\documentclass[11pt,a4paper]{amsart}
\title{Non-K\"ahler manifolds and GIT-quotients}
\date{}
\usepackage{amsfonts, amsopn, amsmath,amssymb}
\usepackage[all]{xy}
\usepackage{ifthen}
\usepackage{color}                    
\usepackage[latin1]{inputenc}

\def\lie{\mathop {\rm Lie}}

\def\re{\mathop {\rm Re}}
\def\im{\mathop {\rm Im}}
\def\smallsetminus{{\scriptstyle \setminus}}

\def\Z{\ensuremath{\mathbb Z}}
\def\C{\ensuremath{\mathbb C}}
\def\P{\ensuremath{\mathbb P}}
\def\R{\ensuremath{\mathbb R}}
\def\Ce{\ensuremath{\C^*}}
\def\Cm{\ensuremath{\C^m}}

\def\Cetm{\ensuremath{(\Ce)^{2m}}}

\def\F{\ensuremath{\mathcal F}}
\def\M{\ensuremath{\mathcal M}}

\def\LVM{\textit{LVM}}
\def\LVMB{\textit{LVMB}}
\def\le{\ensuremath{(\mathcal L, \mathcal E,m,n)}}
\def\lep{\ensuremath{(\mathcal L', \mathcal E',m',n')}}
\def\T{\ensuremath{\mathbb T}}
\def\inv{^{-1}}
\def\cf{{\em cf. }}
\def\ddbar{\partial\bar{\partial}}

\newtheorem{theorem}{Theorem}[section]

\newtheorem{lemma}[theorem]{Lemma}
\newtheorem{proposition}[theorem]{Proposition}
\newtheorem{example}[theorem]{Example}

\newenvironment{remark}
{\trivlist
\item[]\noindent{\bf Remark.}\enspace\ignorespaces}
{\endtrivlist}

\newenvironment{definition}
{\trivlist
\item[]\noindent{\bf Definition.}\enspace\ignorespaces}
{\endtrivlist}

\definecolor{red}{rgb}{1,0,0}
\newboolean{draft}
\setboolean{draft}{false}
\newcommand{\cmt}[1]
{
\ifthenelse{\boolean{draft}}{ {\color{red} #1} }{}
}

\begin{document}
\maketitle
\vspace{-.5cm}
\centerline{\sc St\'ephanie Cupit-Foutou\footnote{\,Universit\"at zu K\" oln, Germany -- scupit@mi.uni-koeln.de --} 
and Dan Zaffran\footnote{\,Academia Sinica, Taipei, Taiwan -- zaffran@math.sinica.edu.tw --}
       }

\vspace{.5cm}
\bigskip

\thispagestyle{empty}

\noindent
L\'opez de Medrano and Verjovsky discovered in 1997 a way to construct many compact
complex manifolds ({\em cf.} \cite{ldm}).
They start with a \C-action on $\P^n$ induced by a diagonal linear vector field
(satisfying certain properties),
and find an open dense subset $U\subset \P^n$ where the action is
free, proper and cocompact, so the quotient $N=U/\C$ is a compact complex manifold.
Their construction was extended to $\C^m$-actions by Meersseman in \cite{m},
yielding a vast family of non-K\"ahler compact manifolds, called \LVM-manifolds.
These manifolds lend themselves very well to various computations,
and a thorough study of their properties is conducted in \cite{m}.
Furthermore, they are (deformations of) a
very natural generalization of Calabi-Eckmann
manifolds. Finally, the topology
of \LVM-manifolds can be extraordinarily complicated: We refer to \cite{bm}
for the most recent results about a study started off in \cite{w} and \cite{ldm}.

It was also remarked that each \LVM-manifold carries
a transversely K\"ahler foliation $\mathcal F$ (\cite{ln}, \cite{m}).

There were two main developments on \LVM-manifolds:
\begin{enumerate}
\item [({\it a\,})]
In \cite{bosio}, Bosio showed that Meersseman's construction
could actually be generalized to more general actions of $\C^m$.
He produced a family of quotient manifolds, containing the family of
\LVM-manifolds,
that we will call \LVMB-manifolds.
He showed that many properties
of \LVM-manifolds carry over to the case of \LVMB-manifolds, and 
also pointed out non-trivial combinatorics relevant to these quotients. 
However, he did not mention
the foliation $\mathcal F$ (that still exists).

\item [({\it b\,})]
In 2004, Meersseman and Verjovsky investigated in \cite{mv}
a link between \LVM-manifolds
and Mumford's Geometric Invariant Theory (Mumford's GIT) that was soon discovered after
\cite{m} was completed. Their main results are that each
\LVM-manifold satisfying {\em condition $(K)$}
(\cf Sect. \ref{LVMBGIT}) admits a Seifert
fibration over a projective simplicial toric
variety $X$, and moreover any such $X$ can appear this way.
\end{enumerate}

\noindent
Building on both ({\it a\,}) and ({\it b\,}), we establish a link between
\LVMB-manifolds and GIT.
We show that the extension from \LVM-manifolds to \LVMB-manifolds
parallels exactly the extension from Mumford's GIT
to the generalized GIT of 
Bia\l ynicki-Birula and \'Swi\c ecicka, 
which allows for non-projective quotients ({\em cf.}  \cite{bbsw2}). 

\smallskip
We describe a construction of \LVMB-manifolds from a GIT point of view,
from which it will be clear that some of them are Seifert-fibered
over a {\em complete} simplicial toric variety $X$. Using a result of
Hamm, we show that any such $X$ can appear that way, i.e., below some \LVMB-manifold.
This generalizes and simplifies the proof of the main result of \cite{mv}.

In most cases $X$ is an algebraic reduction of $N$, and
if $X$ is projective then $N$ is actually an \LVM-manifold.
Using this and the almost-homogeneous structure of these manifolds, 
we produce examples of \LVMB-manifolds that are not  
biholomorphic to any \LVM-manifold. 

\smallskip
Now consider an \LVM-manifold $N$ satisfying condition $(K)$. 
Meersseman and Verjovsky showed that the foliation $\F$ is given 
by a Seifert fibration $N \rightarrow X$, with $X$ projective. 
The fact that $X$ is K\"ahler is the ``reason" for $\mathcal F$ to be transversely K\"ahler.

By our results, an \LVMB-manifold that
is not an \LVM-manifold, and that satisfies condition $(K)$, also admits a
map $N \rightarrow X$, although this time $X$ is {\em not} projective.
We prove that
$\mathcal F$ on such an $N$ is {\em not} transversely K\"ahler. 
The difficulty here is to deal with the singularities
of $X$ as an orbifold.

\smallskip
Finally, our GIT point of view on \LVMB-manifolds leads naturally to a further extension
of the \LVMB\/ family, in the context of the Bia\l ynicki-Birula and Sommese conjecture
(\cf \cite{bbso}). Our results show that the construction of \LVMB-manifolds
actually follows from the solution of this conjecture
for linear algebraic $(\Ce)^k$-actions on the projective space
(\cf \cite{bbsw2}).

\bigbreak
\noindent
\emph{Acknowledgments.}\enspace
A.~Huckleberry invited us as postdocs to Bochum's Ruhr-Universit\"at
where we started this work. We thank him sincerely for his warm hospitality.
We thank also P.~Heinzner, J.~Hubbard, L.~Meersseman and R.~Sjamaar for useful discussions.
We are grateful to P.~Littelmann for having given us the opportunity to
accomplish this work in good conditions.

\section{\LVMB-manifolds}
\label{LVMB}

Let $m$ and $n\geq 2m$ be positive integers. 
Let $\mathcal L = (\ell_0,\ldots,\ell_n)$ be 
an ordered $(n+1)$-tuple of linear forms
on $\C^m$. Fixing an isomorphism between \Cm\ and its dual vector space,
we will look at each $\ell_i$ as a row vector in \Cm.
We will also look at $\ell_i$ as an element in $\mathbb R^{2m}$
via the identification
$\ell_i\mapsto (Re\ \ell_i,Im\ \ell_i) \in \R^m\times \R^m \approx \R^{2m}$.

Let $\mathcal E={\{\mathcal E_{\alpha}\}}_{\alpha}$ be a family of subsets of
$\{0,\ldots, n\}$, each of these subsets
having cardinality $2m+1$. 

We denote these data by $\le$. 
To any given $\le$, we associate the following objects:

-- a Zariski open subset $U=U(\mathcal E)\subset \P^n$ 
defined by 
$$
U=\{[x_0:\dots:x_n]\ \vert\quad 
\text{there exists } \mathcal E_{\alpha}\in {\mathcal E} 
\mbox{ such that for all } i\in \mathcal E_{\alpha} ,\ x_i\neq 0 \};$$

-- a $\C^m$-action on $U$ defined by 
\begin{equation}
\begin{split}
\Cm \times U        &  \longrightarrow  U \\
\left( z,[x_0:\ldots:x_n] \right) &   \longmapsto
[e^{\ell_0(z)}x_0 : \ldots : e^{\ell_n(z)}x_n];
\end{split}\label{hol} \tag{\textsf{hol}}
\end{equation}

-- for all $\alpha$, we define $P_{\alpha}$ as the convex 
hull of $\{ \ell_i,\ i\in \mathcal E_{\alpha} \}$. Notice that 
$2m+1$ is precisely the number of vertices of a simplex in $\R^{2m}$.

\smallskip
\begin{definition}
Let $G$ be a complex Lie group acting holomorphically and effectively on $\mathbb P^n$. 
A $G$-stable open subset $U$ of $\mathbb P^n$ will be called {\it good with respect to
the $G$-action}
when the restricted action on $U$ is proper 
(therefore the topological quotient $U/G$ is Hausdorff), 
and $U/G$ is compact. \hfill$\square$
\end{definition}

\begin{theorem}[\cite{bosio}, 1.4]
\label{bosio_s_theorem}
Take any $\le$, with corresponding $U$ and $\Cm$-action. 
Then $U$ is good with respect to that action 
if and only if the following two conditions are satisfied:
\smallbreak
\noindent
{(\textsf{sep})}\enspace
$\stackrel{\circ}{P_{\alpha}}\cap \stackrel{\circ}{P_{\beta}}  \neq\varnothing$ for every $\alpha, \beta$; 
\smallbreak
\noindent
{(\textsf{comp})}\enspace
for every $\mathcal E_{\alpha}\in \mathcal E$ and every $i\in \{0,\ldots,n\}$
there exists $j\in \mathcal E_{\alpha}$ such that 
$(\mathcal E_{\alpha}- \{j\})\cup \{i\}\in \mathcal E$.
\end{theorem}

When $\le$ satisfies these conditions, 
we write $\le\in\LVMB$, 
and denote the quotient $U/\Cm$ by $N=N\le$. 
By definition, the $\Cm$-action is proper, but 
$\Cm$ has no compact subgroups, so this action is necessarily free.
It is known that a quotient by a free and proper action is a complex manifold
(\cf \cite{Hu} Proposition 2.1.13). Therefore $N$ is a compact complex manifold of dimension $n-m$, 
that 
we call an {\it \LVMB-manifold} (Bosio calls it a \emph{generalized \LVM-manifold}). 

When $\bigcap_{\alpha} \stackrel{\circ}{P_{\alpha}} \neq\varnothing$, the action on $U$
is called an {\it \LVM-action}. We write $\le\in\LVM$, and call $N$ an {\it \LVM-manifold}. 
By Proposition 1.3 in \cite{bosio}, \LVM-manifolds 
are exactly the manifolds constructed in \cite{m}. 

\begin{example}
\label{ex}
{\em Take $n=5$, $m=1$, so the $\ell_i$'s are just complex numbers, 
that --for later use-- we take in $\Z+\sqrt{-1}\;\Z$ according to the diagram: 
$$
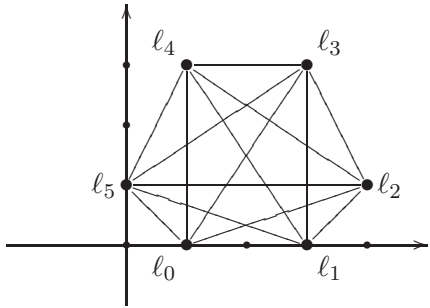
\xy
(-8,0)*{\bullet}="A"+(-3,-3)*{\ell_0}; 
(8,0)*{\bullet}="B"+(3,-3)*{\ell_1}; 
(16,8)*{\bullet}="C"+(3,0)*{\ell_2};
(8,24)*{\bullet}="D"+(3,3)*{\ell_3};
(-8,24)*{\bullet}="E"+(-3,3)*{\ell_4}; 
(-16,8)*{\bullet}="F"+(-3,0)*{\ell_5}; 
"A"; "B" **\dir{-};
"A"; "C" **\dir{-};
"A"; "D" **\dir{-};
"A"; "E" **\dir{-};
"A"; "F" **\dir{-};
"B"; "C" **\dir{-};
"B"; "D" **\dir{-};
"B"; "E" **\dir{-};
"B"; "F" **\dir{-};
"C"; "D" **\dir{-};
"C"; "E" **\dir{-};
"C"; "F" **\dir{-};
"D"; "E" **\dir{-};
"D"; "F" **\dir{-};
"E"; "F" **\dir{-};
(-32,0)*{}="X"; 
(24,0)*{}="Y"; 
(-16,-8)*{}="W"; 
(-16,32)*{}="Z"; 
{\ar@{->} "X"; "Y"};
{\ar@{->} "W"; "Z"};
(-16,0)*{\scriptscriptstyle \bullet};
(-16,16)*{\scriptscriptstyle \bullet};
(-16,24)*{\scriptscriptstyle \bullet};
(0,0)*{\scriptscriptstyle \bullet};
(16,0)*{\scriptscriptstyle \bullet};
\endxy
$$ 
and we take 
$$\mathcal E=\Big\{ 
\{024\}, \{135\},
\{025\},\{124 \}, \{034 \}, 
\{035 \}, \{134 \},\{125 \}   
\Big\}.$$
Then it is easy to check that $\le\in \LVMB$. It is also immediate that $\le\not\in\LVM$
(\cf \cite{bbsw1} p.26 and \cite{bosio} p.1291). 
We will show that the $\LVMB$-manifold $N\le$ is 
not isomorphic to any $\LVM$-manifold.
}
\end{example}

\begin{remark}
{\rm (i)}\enspace
For an action on $\P^n$ of form (\ref{hol}), there are in general many good $U$'s, and on
some of them the restricted action is an \LVM-action, whereas the restricted action on others
is not an \LVM-action. In other words, for $\le\in\LVMB$, there may be $\mathcal E'\neq\mathcal E$ 
such that $(\mathcal L, \mathcal E', m,n)$ is in $\LVMB$ or in $\LVM$.
\smallbreak
\noindent
{\rm (ii)}\enspace
For $\le\in\LVMB$ but $\le\not\in\LVM$, the natural question is whether it is possible that 
$N\le$ be isomorphic to an \LVM-manifold. 
This question was unfortunately left open by \cite{bosio} 
(\cf its review MR1860666-2002i:32015). We will give a partial answer in Sect. \ref{LVMB_not_LVM}. 
\smallbreak
\noindent
{\rm (iii)}\enspace
In the proof of Proposition 2.2 in \cite{bosio}, 
the author shows that if $n> 2m$, $N$ contains a submanifold with odd first
Betti number. From this it follows that $N$ is not K\"ahler.
(Remark: Our $n$ corresponds to $n-1$ in \cite{bosio}.)
In the limiting case
$n=2m$, $N$ is a compact torus, so it is K\"ahler.
\end{remark}

\begin{definition}(\cite{mv})
We say that $\le$ \emph{satisfies condition $(K)$}, and write $\le\in(K)$, when there exists a real
affine automorphism of $(\C^m)^* \approx \C^m \approx \R^{2m}$ sending
each $\ell_i$ to a vector with integer coefficients. 
For example,
an action whose $\ell_i$'s
have only rational coordinates satisfies condition $(K)$. In particular, $(K)$ is a 
dense condition.\hfill$\square$
\end{definition}

\section{Intermediate generalized GIT-quotients}
\label{LVMBGIT}

\noindent
\smallskip
Take integers $m,n$ with $n\geq 2m$. 
%
Consider an algebraic torus $\Cetm$
acting effectively and linearly on $\mathbb P^n$.
Such an action is given by $n+1$
row vectors $\lambda_i = (\lambda_{i,1}, \dots, \lambda_{i,2m})\in \Z^{2m}$ 
for $i=0\dots n$. 
Explicitly:
\begin{equation}
\begin{split}
\Cetm \times \mathbb P^{n}        &  \longrightarrow  \mathbb P^{n}\\
\left( t ,[x_0:\ldots:x_n] \right) &   \longmapsto
[t^{\lambda_0} x_0:\ldots:t^{\lambda_n} x_n]
\end{split}\label{alg}\tag{\textsf{alg}}
\end{equation}
where 
$t^{\lambda_i}$ denotes the product $t_1^{\lambda_{i,1}} \dots t_{2m}^{\lambda_{i,2m}}$.

In this setting, 
Bia\l ynicki-Birula and \'Swi\c ecicka describe 
all the good open subsets in~\cite{bbsw2}. 
This is a generalization of Mumford's GIT, because the 
quotient of such a good open subset $U$ 
is not necessarily projective: 
In general $X=U/(\Ce)^{2m}$ is a complete simplicial toric abstract-algebraic 
variety\footnote{\,These authors also significantly improved Mumford's theory of projective quotients.}. 
In their language, $U\rightarrow X$ is a ``strongly geometric complete quotient" 
(NB: their ``good quotients" are not assumed to be compact). In particular, 
all the \Cetm-isotropy subgroups are finite (\cf~\cite{bbsw2} 1.4).

\smallskip
Now pick $G$ a closed cocompact complex Lie subgroup of
\Cetm\ isomorphic to $\C^m$ (there are plenty of these: see Lemma \ref{closed_cocompact}).
As \Cetm\ acts properly on $U$, so does $G$. Thus $U$ is good with respect to
the $G$-action. 

Moreover, $G$ has no torsion, so it can't intersect the \Cetm-isotropy
subgroups. 
Therefore $G$ acts freely, and $N= U/G$ is a compact complex manifold, 
that we call an {\em intermediate (generalized) GIT-quotient}:
$$
\xymatrix{
U \ar[dd] \ar[rd] & \\
& N=U/G \ar[ld] \\
X=U/\Cetm &
}$$
and the map $N\rightarrow X$ is the quotient
of $N$ by the action of $\T = \Cetm/G$, which is a compact
complex $m$-torus.

As \T\ is compact, the structure of the map $N\rightarrow X$
is very well understood thanks to the results of Holmann
(\cf \cite{o} pp. 82-84). In some cases, \T\ acts freely, so $X$ is
a manifold and $N\rightarrow X$ is simply a \T-principal bundle.
In general, we have a so-called {\em Seifert principal bundle}\,:
$X$ is an orbifold whose singularities correspond to orbits
with non-trivial isotropy. The map $N\rightarrow X$ is ``not
locally trivial'' around those orbits. 
\begin{remark} An intuitive (but {\em real}) 
picture of such an exceptional orbit in a Seifert fibration is provided by 
the circle action $t.(z,w)=(t^2 z, t^3 w)$ 
on the solid torus 
$T=\{ |z|\leq 1 \} \times \{ |w| = 1 \}$. 
This yields a foliation of $T$ for which 
the core of $T$ (i.e., $\{ |z|=0 \} \times \{ |w| = 1 \}$) is 
a fiber which is covered $3$ times by neighboring fibers, which 
are $(2,3)$-torus knots. 
\end{remark}

\begin{theorem}\label{LVMB_K_are_intermediate_GIT_quotients}
Let $\le\in\LVMB\cap (K)$, and $N=N\le$.
Then $N$ can be obtained 
as an intermediate GIT-quotient as above, for some choice 
of a $\Cetm$-action on $\P^n$, a good open
subset $U$ and a subgroup $G$. In particular, 
any \LVMB-manifold can be obtained as a small
deformation of such an intermediate GIT-quotient.
\end{theorem}

\begin{lemma}
\label{closed_cocompact}
Let $G\approx \Cm$ be a Lie subgroup of \Cetm. Let $A\in \C^{2m\times m}$ be a matrix
whose columns $A_1,\dots, A_m$ form a \C-basis of\, $\lie(G)$. Then:\\
{\rm (i)}\enspace If $G$ is closed then $G$ is cocompact;\\
{\rm (ii)}\enspace $G$ is closed if and only if the matrix
$(\re\ A|\im\ A) \in \R^{2m\times 2m}$
is invertible.
\end{lemma}

\begin{proof}
Denote the exponential map $exp: \C^{2m}\rightarrow\Cetm$, whose
kernel is $i\Z^{2m}$. As $G\approx \Cm$, $\lie(G) \cap i\Z^{2m}= \{0\}$.
Then $G$ is closed if and only if $\lie(G) \cap i\R^{2m}= \{0\}$.
Thus, when $G$ is closed, $\mathrm{Span}_{\Z}(\alpha \cup \beta)$ is a full lattice in
$\C^{2m}$ for any real bases $\alpha$ and $\beta$ of $\lie(G)$ and $i\R^{2m}$
respectively. Therefore a closed $G$ is cocompact.

Now:
$\lie(G) \cap i\R^{2m} \neq \{0\}$ if and only if
there exists $(z_1,\dots,z_m) \in \Cm\smallsetminus \{0\}$ such that
\begin{equation}
\label{inversibility}
\re\ \Big(\sum_{j=1}^m z_j A_j\Big)=0.\tag {$\star$}
\end{equation}
Writing $z_j$ as $x_j+iy_j$, the equation~(\ref{inversibility}) is equivalent to
$$\sum_{j=1}^m x_j \re\ A_j -y_j \im\ A_j=0.$$
Thus $G$ is not closed if and only if $\{\re\ A_1, \dots, \re\ A_m, \im\ A_1, \dots, \im\ A_m \}$ is \R-linearly dependent,
that is, $(\re\ A|\im\ A)$ is not invertible.
\end{proof}

\begin{proposition}
\label{factorisation}
An action of form~(\ref{hol})
satisfies $(K)$ if and only if
it is the restriction of a \Cetm -action of
form~(\ref{alg})
to a closed cocompact subgroup isomorphic to $\mathbb C^m$.
\end{proposition}

\begin{proof}
Take a \Cetm-action of form (\ref{alg}) given
by $\lambda_0,\dots, \lambda_n \in \Z^{2m}$.
Let $G$ be a closed cocompact subgroup of \Cetm\/
which is isomorphic to \Cm.
Take
$A \in \C^{2m\times m}$ as in Lemma \ref{closed_cocompact}.
We know that $(\re\ A|\im\ A)$ is invertible.
Writing $G=\{ \exp(Az),\ z\in\Cm \}$, a direct
computation shows that the \Cetm-action restricted to $G$ is
of form (\ref{hol}), with $\ell_i= \lambda_i A$ for $i=0\dots n$.
That restricted action satisfies $(K)$ because
$ (\re\ \ell_i,\im\ \ell_i) (\re\ A|\im\ A)\inv
\in \R^m \times \R^m\approx \R^{2m}$
equals $\lambda_i$, so has integer coefficients.

Conversely, take a  \Cm-action of form (\ref{hol}) satisfying $(K)$,
given by the linear forms $\ell_0, \dots, \ell_n$. Notice that translating all the
$\ell_i$'s by the same vector does not change the action.
There exist $\lambda_0,\dots, \lambda_n \in \Z^{2m}$, an invertible
$M \in \R^{2m\times 2m}$ and $b \in \R^{2m}$ such that for all $i$,
$(Re\ \ell_i, Im\ \ell_i) = \lambda_i M +b$. Now take $A\in \C^{2m\times m}$
such that
$(Re\ A|Im\ A)= M$, and denote by $G$ the subgroup of \Cetm\ such that
$\lie(G)= \mathrm{Image}(A)$. By Lemma \ref{closed_cocompact},
$G$ is closed and cocompact. Moreover the
restriction to $G$ of the \Cetm-action given by the $\lambda_i$'s is
the \Cm-action we started with.
\end{proof}

\noindent
{\it Proof of Theorem \ref{LVMB_K_are_intermediate_GIT_quotients}.}
Let $\le\in\LVMB\cap (K)$, with associated $N=U/\mathbb C^m$.
Take a $\Cetm$-action
of form (\ref{alg}) given by Proposition \ref{factorisation}.

By definition, $U$ is the complement of some coordinate subspaces, so
it is stable by the ``big'' algebraic torus $(\Ce)^{n+1}$ acting diagonally on $\P^n$,
so in particular $U$ is stable by that \Cetm-action.
Moreover 
the topological quotient 
space $U/\Cetm$ is homeomorphic to $X=\big(U/\C^m\big)\big/ \big(\Cetm/\C^m\big)$,
which
is the quotient of the compact manifold $N$ by a compact group, 
therefore it is compact and Hausdorff. Hence $U\rightarrow X$ 
is one of the ``strongly geometric complete" quotients 
considered in~\cite{bbsw2}. 
\hfill$\square$

\section{Applications}
\label{applications}\noindent
{\bf Notation.} For a given $\le$, we denote by $d$ the minimal codimension of $\P^n - U$ 
(\cf \cite{m}). The condition $d>1$ is equivalent to $\bigcap_{\alpha} \mathcal E_{\alpha} = \varnothing$.
\hfill$\square$

\subsection{\LVMB-manifolds and complete toric varieties}
\label{LVMB-manifolds_and_complete_toric_varieties}
An action of form (\ref{hol}) (resp. (\ref{alg})) determines 
a subgroup $H_1$ (resp. $H_2$) of $(\Ce)^{n+1}$, 
with $H_1\approx \Cm$ and $H_2\approx (\Ce)^{2m}$. 
When $(K)$ is satisfied, 
the extension from a $\Cm$-action to a $(\Ce)^{2m}$-action 
in Proposition \ref{factorisation}
is unique in the sense that
$H_2$ must be 
the smallest algebraic subtorus of $(\Ce)^{n+1}$ containing $H_1$, 
as follows 
from the cocompactness of $H_1$ in $H_2$. 
In particular, for any given $\le\in\LVMB \cap (K)$, 
there is a well-defined $$X=X\le := U/H_2,$$ with a map $\pi:N\le \rightarrow X$ 
which is a Seifert principal fibration with compact torus $\T=H_2/H_1$, 
called a {\em generalized Calabi-Eckmann fibration} in \cite{mv}.
\cmt{check coherence} 
The fibers of $\pi$ define a foliation on $N$ denoted by $\F=\F\le$. 

The next result extends one of \cite{mv}'s main results to a more general setting. We also
give a simpler proof, based on a result of H. Hamm.

\begin{theorem}
\label{calabi}
Let $X$ be any complete simplicial
toric variety. Then there exists an \LVMB-manifold $N$
giving a generalized
Calabi-Eckmann fibration over $X$.
\end{theorem}

\begin{proof}
Applying Theorem~6.1 in~\cite{hamm},
one can realize $X$
as a geometric quotient $V/(\Ce)^r$
with $V$ an open subset of some $\mathbb C^n$
acted on linearly by $(\Ce)^r$ and $n>r$. This action
is the restriction (to a subgroup isomorphic to $(\Ce)^r$)
of the diagonal  $(\Ce)^n$-action on $\C^n$. If $r$ is
odd, we let $(\Ce)^{r+1}= (\Ce)^r \times \Ce$ act on
$V \times \Ce$ by $(t_1, \dots, t_r, t).(v, z)= ((t_1, \dots, t_r).v, tz)$
then $V \times \Ce$ is an open subset of $\C^{n+1}$, and
$V\times \Ce /(\Ce)^{r+1} = X$.
Therefore (up to replacing $V$ with $V \times \Ce$ and $n$ with $n+1$) we can
assume that $r$ is even, i.e., $r=2m$.

The map 
$(z_1,\dots, z_n)\mapsto [1:z_1:\dots:z_n]$ defines a 
\Cetm-equivariant embedding $\C^{n} \hookrightarrow \P^n$ sending 
$V$ to an open
subset $U$ of $\P^n$ whose quotient by $\Cetm$ is $X$. Taking any 
$G\approx \Cm$ closed and cocompact in \Cetm\ gives a generalized
Calabi-Eckmann fibration $N=U/\Cm$ above $X$.
\end{proof}

\begin{proposition}
\label{NLVMssiXproj}
Let $\le\in\LVMB\cap (K)$, with corresponding $N$ and $X$. 
\smallbreak
\noindent
{\rm (i)} 
If $d>1$ then $X$ is the algebraic reduction of $N$;
\smallbreak
\noindent
{\rm (ii)} 
$X$ is projective if and only if $\le\in\LVM$.\\ 
\phantom{{\rm (ii)}}
In particular, if $X$ is projective then $N$ is an \LVM-manifold.
\end{proposition}

\begin{proof}

(i) The proof of Theorem~4 in~\cite{m} applies here
and shows that any meromorphic function on $N$ is 
a function of some functions $M_1,\dots, M_a$. 
But $M_1,\dots, M_a$ are pull-backs
of meromorphic functions on $X$ by Remark 2.13 in \cite{mv}.
As $X$ is an abstract-algebraic variety,
it is a (non-necessarily projective) algebraic reduction of $N$.

(ii) Take any \Cetm-action of form (\ref{alg}) given by 
Proposition \ref{factorisation}. 
From the proof of Proposition \ref{factorisation}, we know 
that the $\lambda_i$'s of this action are obtained from the 
$\ell_i$'s of the \Cm-action by an affine automorphism, that we denote $\varphi$. 
For each polytope $P_{\alpha}$ associated to $\le$, we define $Q_{\alpha}=\varphi(P_{\alpha})$. 

Denote $\Pi= {\{ Q_{\alpha} \}}_{\alpha}$, and define $\tilde{U} \subset \P^n$ by: 
$[x_0:\dots:x_n]\in \tilde{U}$ if and only if the convex hull of 
$\{\lambda_i\ \vert\ x_i\neq 0\}$ 
contains a polytope of $\Pi$. 
We claim that $\tilde{U}=U$. 
The inclusion $U\subset \tilde{U}$ follows from the definitions. 
On the other hand, 
it follows from Theorem \ref{bosio_s_theorem} (\textsf{sep}) that 
\begin{equation}
\label{interiors_intersect}
\text{for all } Q,Q' \in\Pi, \quad 
\stackrel{\circ}{Q} \cap \stackrel{\circ}{Q'} \neq \varnothing.
\tag {$\cap$}
\end{equation}
In particular, any polytope of $\Pi$ has dimension $2m$, so the quotient 
of $\tilde{U}$ by $\Cetm$ is a geometric quotient (\cf Theorem~7.8 in \cite{bbsw2}), i.e., 
the orbit space $\tilde{U}/\Cetm$ is Hausdorff. But this space also contains the compact set $U/\Cetm$ as 
an {\em open} subset, so $\tilde{U}/\Cetm = U/\Cetm$. 
This proves the claim, so $U$ is given by $\Pi$ in the sense of \cite{bbsw2}. 

From Example 7.12 and Corollary 7.16 in~\cite{bbsw2}, 
it follows that a necessary and sufficient condition for 
$X$ to be projective is that $\bigcap_{Q\in\Pi} Q \neq \varnothing$. 
By (\ref{interiors_intersect}) this is 
equivalent to 
$$\bigcap_{Q\in\Pi} \stackrel{\circ}{Q} \neq \varnothing,$$
which is in turn equivalent to $\bigcap_{\alpha} \stackrel{\circ}{P}_{\alpha} \neq\varnothing$,
i.e., $\le\in\LVM$. 
\end{proof}

\subsection{Existence of \LVMB-manifolds that are not \LVM-manifolds}
\label{LVMB_not_LVM}

\subsubsection{Standard submanifolds}
Let $\le\in\LVMB$, with associated $N=N\le$. 
Let $U^1$ be a subset of $U\subset\P^n$ defined by the vanishing of some 
homogeneous coordinates. 
Then the $\Cm$-action on $U$ restricts on $U^1$ to an action of form (\ref{hol}), 
and when it satisfies the conditions of Theorem \ref{bosio_s_theorem}, it gives 
a submanifold $N^1\subset N$ called 
a {\em standard submanifold
of $N$ with respect to $\le$} (\cf \cite{m} p.100).

Examples of standard submanifolds
are obtained as follows:
Take
any element $\mathcal E_{\alpha}$ of $\mathcal E$, and define a $2m$-dimensional 
$U^1$ by letting 
$z_i=0$ for $i\not\in \mathcal E_{\alpha}$. 
By doing this for each
element of $\mathcal E$, one gets
a finite family of standard $m$-submanifolds. They have
minimal dimension among standard submanifolds of $N$ w.r.t. $\le$.
We call them {\em minimal} standard submanifolds of $N$ w.r.t. $\le$.

\begin{lemma} \label{minimal_iff_torus}
A standard submanifold of $N$ 
is minimal
if and only if it is a compact complex torus.
\end{lemma}
\begin{proof}
Each standard submanifold is itself an \LVMB-manifold,
so the statement follows from Proposition 2.1 in \cite{bosio},
and the fact that tori are symplectic (NB: our $n$ corresponds to $n-1$ in \cite{bosio}).
\end{proof}

\subsubsection{Action of $Aut^0(N)$}\label{action_of_aut0}
We state a few consequences of \cite{bosio} Proposition 2.4 and of \cite{m} Part V, 
whose proofs apply also to \LVMB-manifolds. 

Let $N=N\le$ be an $\LVMB$-manifold. 
The linear diagonal action of $(\Ce)^{n+1}$ on $\P^n$ descends
to an action on $N$, which corresponds to a subgroup of 
$Aut(N)$ that we call $G$. Each $G$-orbit is a standard submanifold $N^1$
minus all standard submanifolds strictly contained in $N^1$.
In particular, minimal standard submanifolds are $G$-orbits.

In general, $G$ is contained in $Aut^0(N)$, the connected component
of the identity. Therefore the closure of any $Aut^0(N)$-orbit is
a union of standard submanifolds.

Finally, if $d>1$ and $\ell_i\neq\ell_j$ when $i\neq j$,
then $G=Aut^0(N)$. This follows from \cite{m} Theorem 6, which states the
result at the --commutative-- Lie algebra level.

\subsubsection{A family of \LVMB-manifolds that are not \LVM-manifolds}

\begin{proposition}\label{mnK_intrinsic}
Let $\le\in\LVMB$, with associated $d$
and $\ell_i$'s such that $d>1$ and $\ell_i\neq \ell_j$ when $i\neq j$. 
Let $N=N\le$. 

If there exists $\lep$ such that $N\lep = N$, then:
\smallbreak
\noindent
{\rm(i)}\enspace $m'=m$, $n'=n$, $d'>1$;
\smallbreak
\noindent
{\rm(ii)}\enspace if $\le\in (K)$ then $\lep\in (K)$.
\end{proposition}

\begin{proof}
We use the facts stated in \ref{action_of_aut0}.

(i) The $Aut^0(N)$-orbits of
minimal dimension are exactly the minimal standard submanifolds
of $N$ w.r.t. $\le$. Therefore they are tori of dimension $m$.

Let $p$ belong to one of these tori. We know {\em a priori } that 
the closure of the orbit
of $p$ under $Aut^0(N)$ is a union of standard submanifolds
w.r.t. $\lep$. That orbit being closed and smooth, 
it must equal exactly {\em one} such
submanifold $N'_i$. By Lemma \ref{minimal_iff_torus}, it is a
minimal standard submanifold w.r.t. $\lep$, so $m'=m$.
Then $n'=n$ because $n'-m'=dim\ N= n-m$.

From the proof of Proposition 2.1 in \cite{bosio}, $d>1$ implies that 
$$H^2(N,\Z)\stackrel{(*)}{=}\Z \text{ and } \pi_1(N)=0$$ ($\pi_1(N)=0$ because there is
a principal bundle over $N$ with fiber a circle and a
$2$-connected total space).

Now we use that $N=N\lep$.
Denote by $k'$ the cardinality of $\bigcap_{\alpha} \mathcal E'_{\alpha}$.
Assume that $d'=1$. Then $k'\geq 1$, and Bosio's description
of the homotopy type of $N$ implies that $H^2(N,\Z)$ has rank $0$ or $(k'-1)(k'-2)/2$
and $\pi_1(N)=\Z^{k'-1}$. Thus, by $(*)$, $(k'-1)(k'-2)/2=1$, so $k'=3$.
Then $\pi_1(N)=\Z^2$, which is a contradiction, so $d'>1$.

(ii) For $N=N\le$ with $d>1$, Theorem 4 (ii) in \cite{m}
equates the algebraic dimension
of $N$ to an integer $a$ (the proof given for \LVM-manifolds
applies also in the \LVMB\/ case). It follows from the definitions of $a$ 
and $(K)$ that in general 
$a\leq \dim N - m$, with equality if and only if $\le\in (K)$.
Therefore (ii) follows from (i).
\end{proof}

\begin{proposition}\label{F_intrinsic}
Let $N=N\le$ be an \LVMB-manifold. 
Assume that $\le\in (K)$, $d>1$,
and $\ell_i\neq \ell_j$ for $i\neq j$.

Then there exists an open, dense and $\F$-saturated subset $N_* \subset N$
such that the following are equivalent:
\smallbreak
\noindent
{\rm(i)}\enspace $p,q \in N_*$ belong to the same leaf of $\F$.
\smallbreak
\noindent
{\rm(ii)}\enspace for all $f\in \M(N)$ holomorphic at $p$ and $q$, $f(p)=f(q)$.
\end{proposition}

\begin{proof}
The foliation $\F$ is given by the fibers of the map
$\pi: N\rightarrow X\le$ and $X$ is an abstract-algebraic variety,
so in particular $X$ is a Moishezon space. By Hironaka's and Moishezon's 
theorems, there exists
$\tilde{X}$ smooth and projective,
with a sequence of blow-downs $\varphi: \tilde{X} \rightarrow X$.
$$\xymatrix{
                & \tilde{X} \ar[d]_{\varphi}   & \tilde{X}_* \ar@{_{(}->}[l] \ar[d]^{\sim}\\
N \ar[r]^{\pi} & X                            & X_*         \ar@{_{(}->}[l]\\
N_* \ar@{^{(}->}[u] \ar[rru]
           }$$
Then $f\mapsto f\varphi$ is an isomorphism between the fields
of meromorphic functions $\mathcal M(X)$ and $\mathcal M(\tilde{X})$,
and there exist $X_*$ and $\tilde{X}_*$ open dense subsets in
$X$ and $\tilde{X}$ between which $\varphi$ induces a biholomorphism.
Define $N_*=\pi^{-1}(X_*)$.

(ii) implies (i):\\
Take $p,q\in N_*$ not on the same leaf, i.e.,
$\pi(p)$ and $\pi(q)$ are distinct points in $X_*$,
so $p_*= \varphi^{-1}\pi(p)$ and $q_*= \varphi^{-1}\pi(q)$
are distinct in $\tilde{X}_*$.
Then there exists $f\in\mathcal M(\tilde{X})$ holomorphic
at $p_*$ and $q_*$ such that $f(p_*)\neq f(q_*)$ (\cf \cite{gf}, V, Theorem 3.14).
Then $f\varphi^{-1}\pi\in \mathcal M(N)$ is holomorphic
at $p$ and $q$, and $f(p)\neq f(q)$.

(i) implies (ii):\\
Take $p,q\in N_*$ on the same leaf $L$, i.e.,
$\pi(p) = \pi(q)$. Let $f\in \mathcal M(N)$, holomorphic at $p$ and $q$. 
Then there exists $f_1\in \mathcal M(X)$ such that $f=f_1\pi$
(\cf proof of Proposition \ref{NLVMssiXproj}).

As $f$ is holomorphic at $p$, it is bounded
on a neighborhood of $p$. As $\pi$ is an open map,
$f_1$ is bounded on a neighborhood $A$ of $\pi(p)=\pi(q)$.
So $f_1$ is holomorphic on $A$, and $f$ is holomorphic on
$\pi^{-1}(A)$, and in particular on $L$, which is compact.
So $f$ is constant on $L$, and in particular $f(p)=f(q)$.
\end{proof}

\begin{theorem}\label{LVMBs_not_LVM}
Let $\le\in \LVMB \cap (K)$, with $d>1$ and $\ell_i\neq \ell_j$ for $i\neq j$.
If $\le\not\in \LVM$ then $N$ is not biholomorphic to any \LVM-manifold.
\end{theorem}
\begin{proof}
Assume on the contrary
that $N$ can also be written as $N\lep$ with $\lep\in\LVM$.
By Proposition \ref{mnK_intrinsic} (ii), $\lep\in(K)$.
Now denote by $\F$ and $\F'$ the foliations of $N$ with respect to $\le$ and $\lep$.
By Proposition \ref{F_intrinsic}, $\F$ and $\F'$ agree on a dense open subset,
so $\F=\F'$ everywhere. But by Proposition \ref{NLVMssiXproj},
$X'=N/\F'$ is projective, whereas $X=N/\F$\, is not. Contradiction.
\end{proof}

\noindent
{\bf Example.} In Ex. \ref{ex}, it is immediate to check that $\le\in(K)$ and that 
$\bigcap_{\alpha} \mathcal E_{\alpha} = \varnothing$, so $d>1$ (actually $d=2$). 
Also $\le\not\in\LVM$, so
by Theorem \ref{LVMBs_not_LVM}, $N\le$ is not isomorphic to any $\LVM$-manifold.

\subsection{Foliations}
\label{foliations}

Let $\le\in\LVM$, and $N=N\le$. 
Generalizing the results of L\oe b and Nicolau, Meersseman shows ({\em cf.} \cite{m})
the existence of a foliation $\F=\F\le$ on $N$, and proves it is transversely
K\"ahler. This means that $N$ admits a closed, real and $J$-invariant
two-form $\omega$, positive on the normal bundle of $\F$, and such that $\ker(\omega)=\F$.
This is a strong property, that has interesting geometric consequences on $N$. 

\smallskip
Now assume that $\le\in(K)$. Then the foliation's leaves are just 
the fibers of $\pi: N \rightarrow X\le$. 
Assume moreover that $\le\not\in \LVM$. 
By Proposition \ref{NLVMssiXproj} (ii), we know
that $X$ is not projective. Using this fact, we prove below 
that $\F$ is not transversely K\"ahler. 
This is an unexpected difference
with the case of \LVM-manifolds. 

The proof is straightforward when $\T$ acts freely on $N$, i.e., 
when $N\rightarrow X=N/\T$
is a genuine principal bundle: Assume that \F\ is transversely
K\"ahler with respect to a 2-form $\omega$. Then:
\begin{enumerate}
\item[(a)] Make $\omega$ a \T-invariant form by averaging it
over the \T-action, and push it forward on $X$, where
it gives a K\"ahler form.
\item[(b)] As $X$ is smooth, K\"ahler and Moishezon, 
it is projective by a theorem of B. Moishezon.
\end{enumerate}

In the general case of $N\rightarrow X$ being a Seifert
principal bundle (and $X$ being singular), both steps
(a) and (b) become non-trivial:

For (a): The problem is that the push forward
of a \T-invariant $\omega$ is not smooth in general (it is not
even continuous). To fix this problem, we
use some results and methods of D. Barlet and J. Varouchas. We take
local potentials of $\omega$ on slices of
the Seifert bundle, and push them forward to $X$.
Then we can apply the following theorem of Varouchas:
{\em If a complex space $X$ has an open cover $\{U_i \}_{i\in I}$,
with for all $i$ a continuous stricly plurisubharmonic
function $\psi_i$, such that for all $i$ and $j$, $\psi_i - \psi_j$
is pluriharmonic, then $X$ is a K\"ahler space in the sense
of Grauert}. (Grauert's definition of a K\"ahler
space is exactly the above sentence with ``continuous''
replaced with ``smooth''.)

For (b): It is not true
in general that a Moishezon K\"ahler space is
projective. However, as $X$ has only rational singularities,
a theorem of Y. Namikawa implies that $X$ is projective.

\begin{theorem}
Let $\le\in\LVMB\cap (K)$, with corresponding $N$ and $\F$. 
If $\le\not\in\LVM$ then \F\
is not transversely K\"ahler. 
\end{theorem}

\begin{proof}
Assume on the contrary that there exists a 2-form $\omega_{0}$ with respect to which
\F\ is transversely K\"ahler. We have to show that $\le\in\LVM$. 

\smallskip
{\em First step: ``make $\omega_0$ \T-invariant''.}
Let $A$ be an automorphism of $N$ induced by some element
of \T. Remark that for any $p\in N$ and any vector $v\in T_p N$,
$v$ is tangent to the \T-orbit of $p$ if and only if $A_*v$ is.

For $u,v \in T_pN$, define
$$
\omega(u,v)= \int_{t\in \T} \omega_{0}(tu,tv)
$$
for the normalized Haar measure on \T.
This defines a \T-invariant, closed, real and $J$-invariant
2-form $\omega$. Moreover, by the above remark,
$Ker\ \omega=\F$ and $\omega$ is positive on the
normal bundle of $\F$ (because an integral of positive
numbers is positive). So \F\ is transversely K\"ahler
with respect to $\omega$.

\smallskip
{\em Second step: ``push forward local potentials on slices''.}
It follows from the results of Holmann (\cf \cite{o} pp. 82--84)
that we can find a family
of local holomorphic
slices $\{S_i \}_{i\in I}$ such that for each $i$:
\begin{itemize}
\item $S_i \subset N$ is transverse to \F\ and biholomorphic to a ball
of same dimension as $X$;
\item $\pi_i = \pi_{|S_i}$ is a quotient by a finite subgroup
of \T\ denoted by $\Gamma_i$ ($\Gamma_i$ is some isotropy subgroup);
\item $\{V_i = \pi_i(S_i) \}_{i\in I}$ is an open cover of $X$.
\end{itemize}
For each $i$, $\pi_i : S_i \rightarrow V_i$ is a ramified covering.
We denote by $\pi_i: S_i^{reg} \rightarrow V_i^{reg}$ the
associated regular covering (off the ramification locus).
We define $\omega_i = \omega_{|S_i}$, which is a K\"ahler
form on the ball $S_i$. So it admits a potential function
$\varphi_i$, i.e., $\varphi_i \in C^{\infty}(S_i, \mathbb{R})$, $\varphi_i$ is strictly
plurisubharmonic (p.s.h.)
and $\omega_i= \sqrt{-1} \ddbar\varphi_i$. Now define for all $i$
a map $\psi_i : V_i \rightarrow \mathbb{R}$ by
$$
\psi_i(x)= \frac{1}{n_i}\sum_{\gamma \in \Gamma_i} \varphi_i(\gamma . p),
$$
where $p\in \pi_i\inv(x)$ and $n_i$ is the order of $\Gamma_i$.
Then for all $i$, $\psi$ is continuous and strictly p.s.h.
by a result of Barlet ({\em cf.} \cite{v2} Proposition 3.4.1).

\smallskip
{\em Third step: we prove that for all $i,j$,
$\psi_i - \psi_j$ is pluriharmonic on $V_i \cap V_j$.}
We use the definition of pluriharmonic (p.h.) of
\cite{v2}. For a real function, this means that the function is locally the real part of a holomorhic
function.

We will first prove that $\psi_i - \psi_j$ is p.h. on
$V_i^{reg} \cap V_j^{reg}$. Let $V$ be any small open
subset of $V_i^{reg} \cap V_j^{reg}$ isomorphic to an open ball.
Then $\pi_i\inv(V)$ is a union of balls $B_{i,1}, \dots, B_{i,n_i}$
and $\pi_j\inv(V)$ is a union of balls $B_{j,1}, \dots, B_{j,n_j}$.
\def\Ba{\ensuremath{B_{\alpha}}}
\def\Bb{\ensuremath{B_{\beta}}}
\def\ua{\ensuremath{u_{\alpha}}}
\def\ub{\ensuremath{u_{\beta}}}
\def\va{\ensuremath{v_{\alpha}}}
\def\vb{\ensuremath{v_{\beta}}}
\def\pa{\ensuremath{p_{\alpha}}}
\def\pb{\ensuremath{p_{\beta}}}
\def\pia{\ensuremath{\pi_{\alpha}}}
\def\pib{\ensuremath{\pi_{\beta}}}
Pick any two of them \Ba\ and \Bb. Denote $\pia = \pi_{|\Ba}$
$\pib = \pi_{|\Bb}$. Then \pia\ (resp. \pib) sends
\Ba\ (resp. \Bb) isomorphically onto $V$.

We now check that
\begin{equation}
(\pia)_*\omega = (\pib)_*\omega\tag{$\diamond$}.
\end{equation}
Take $p\in V$ and $u,v \in T_pV$.
Denote: $\pa= \pia\inv(p)$, $\ua = (\pia\inv)_*u$, $\va= (\pia\inv)_*v$,
and \pb, \ub, \vb\ in the analogous way. We need to show:
\begin{equation}
\omega(\ua,\va)= \omega(\ub, \vb)\label{one}.
\end{equation}
As \pa\ and \pb\ are on the same fiber of $\pi$, there exists
$A\in Aut(N)$ induced by some $t\in \mathbb T$
such that $A(\pa)=\pb$. By the $\mathbb T$-invariance of $\omega$,
\begin{equation}
\omega(\ua,\va)= \omega(A_*\ua, A_*\va)\label{two}.
\end{equation}
On the other hand, $\pi A=\pi$,
so $\pi_*A_*\ua=\pi_*\ua=u=\pi_*\ub$. Therefore $\ub-A_*\ua \in Ker\ \pi_*=Ker\ \omega$.
Thus $\omega(\ub-A_*\ua, \vb)=0$,
i.e., $\omega(A_*\ua, \vb)=\omega(\ub, \vb)$.
Repeating this for the second entry we get
\begin{equation}
\omega(A_*\ua, A_*\va)=\omega(\ub, \vb) \label{three}.
\end{equation}
By (\ref{two}) and (\ref{three}), we get (\ref{one}). Therefore ($\diamond$)
holds,
which means that pushing forward $\omega$ by $\pi$ from any ball among
$B_{i,1}, \dots, B_{i,n_i}, B_{j,1}, \dots, B_{j,n_j}$ gives on $V$ the same 2-form,
that we denote by $\omega_V$.

By the definition of $\psi_i$ and $\psi_j$ we have
$$
\psi_{i_{|V}}= \frac{1}{n_i}\sum_{\nu= 1\dots n_i} \pi_*(\varphi_{i_{|B_{i,\nu}}})
\ \text{ and }\ 
\psi_{j_{|V}}= \frac{1}{n_j}\sum_{\nu= 1\dots n_j} \pi_*(\varphi_{i_{|B_{j,\nu}}}).
$$
So
\begin{align*}
\sqrt{-1}\ddbar\psi_{i_{|V}}&
=\frac{1}{n_i}\sum_{\nu= 1\dots n_i} \sqrt{-1}\ddbar \pi_*(\varphi_{i_{|B_{i,\nu}}})\\
&=\frac{1}{n_i}\sum_{\nu= 1\dots n_i} \pi_*(\sqrt{-1}\ddbar \varphi_{i_{|B_{i,\nu}}})
\ \ \text{ (because $\pi$ is holomorphic)}\\
&=\frac{1}{n_i}\sum_{\nu= 1\dots n_i} \pi_*(\omega_{|B_{i,\nu}})\\
&= \omega_V,
\end{align*}
and, similarly, $\sqrt{-1}\ddbar\psi_{j_{|V}}= \omega_V$.
Therefore $\ddbar(\psi_i - \psi_j)_{|V}= 0$, so $\psi_i - \psi_j$ is on $V$ the real part
of some holomorphic function.
This proves that  $\psi_i - \psi_j$ is p.h. on $V_i^{reg} \cap V_j^{reg}$,
and we also know that it is continuous on $V_i \cap V_j$.

\smallskip

Now take any small open subset $V\subset V_i \cap V_j$ such that $\pi\inv(V)$
has a simply connected component, that we denote by $U$.
Denote by $Q$ the singular locus of $V$ as an orbifold, and
denote by $R=\pi_i\inv(Q)\cap U$. Then $R$ is the ramification locus of $\pi_i:U\rightarrow V$.
In particular, $R$ is an analytic subset of $U$.
Now $\pi_i^*(\psi_i - \psi_{j_{|V}})$ is continuous on $U$, and it is p.h. on $U-R$, so
in particular it is p.s.h. on $U-R$. By a result of Grauert and Remmert (in \cite{gr}),
that function extends as a p.s.h. function on $U$. We can show similarly that its
opposite is p.s.h. on $U$. Therefore it is p.h., so it is the real part
of some holomorphic function. Now pushing forward that
holomorphic function gives on
$V$ a holomorphic function whose real part is $\psi_i - \psi_j$.

\smallskip
{\em Fourth step:} We can now apply Theorem 1 of \cite{v2}
(remark that our $X$ is reduced, so ``(ii) follows from (i)''),
to get that $X$ is K\"ahler in the sense of Grauert ({\em cf.} \cite{g}). As $X$ is a toric
variety, we know it has only rational singularities by Theorem 5.2 in \cite{c}.
Also $X$ is an abstract-algebraic variety, 
so it is a Moishezon complex space. By Corollary 1.7 in \cite{n},
we get that $X$ is projective. By Proposition \ref{NLVMssiXproj}, $\le\in\LVM$.
\end{proof}

\begin{remark}
For an \LVMB-manifold obtained by a non-\LVM action that does not
satisfy $(K)$, one can still construct a foliation \F, as in \cite{m}.
We expect \F\ to be non transversely K\"ahler, but are not able
to prove it. There is no general deformation
argument to get this result: Whereas a small deformation of a compact
K\"ahler manifold is still K\"ahler, the analogous statement is
not true in general for transversely K\"ahler foliations, unless the differentiable
type of the foliation is preserved ({\em cf.} \cite{ekag}).
\end{remark}

\subsection{Further generalizations of the \LVMB\/ family}

As a concluding remark, we give two ways of generalizing the family of \LVMB-manifolds
that are very natural from our GIT point of view. These two ways are
independent, and can be combined.

\subsubsection{Other cocompact subgroups}
Start with a $(\Ce)^k$-action of form (\ref{alg}),
with $k$ not necessarily even. Now take {\em any}
closed cocompact complex Lie subgroup $G\subset (\Ce)^k$. Such
a $G$ is isomorphic to $\Cm \times (\Z)^l$, with $k=2m+l$.
Now take a good open subset $U$ for the $(\Ce)^k$-action.
The quotient $N=U/G$ is a compact complex manifold. One
can still deform the parameters of the action to get
other manifolds.

Note that $N$ is topologically a fiber bundle
over an \LVMB-manifold with fiber a real torus $(S^1)^l$.

\subsubsection{Quotients \`a la Bia\l ynicki-Birula and Sommese}

Bia\l ynicki-Birula and Sommese have conjectured that part of
the algebraic theory of GIT-quotients can be extended to
the more general case of meromorphic actions of $(\Ce)^k$
on a reduced compact normal complex analytic space $Y$, with similar
combinatorial properties
(\cf \cite{bbso}). Unfortunately, only the cases of
$k=1$, and $k=2$ with $Y$ smooth and K\"ahler are fully
understood so far, and there are few non-toric examples
worked out in the literature.

It is likely that these quotients
yield many new generalized Calabi-Eckmann fibrations, possibly
giving new examples over simplicial toric varieties as well.

\end{document}